\documentclass{amsart}

\newtheorem{theorem}{Theorem}[section]

\theoremstyle{definition}

\newtheorem{example}[theorem]{Example}

\theoremstyle{remark}

\numberwithin{equation}{section}
\DeclareMathOperator{\codim}{codim}

\begin{document}

\title{Classification of rationally elliptic toric orbifolds}

%    Remove any unused author tags.

%    author one information
\author{Michael Wiemeler}
\address{Mathematisches Institut\\WWU M\"unster\\Einsteinstr. 62\\D-48149 M\"unster\\Germany}
%\curraddr{}
\email{wiemelerm@uni-muenster.de}
\thanks{The research for this paper was supported by SFB 878 Groups, Geometry and Actions and the Cluster of Excellence Mathematics M\"unster at WWU M\"unster.
}

%    author two information
%\author{}
%\address{}
%\curraddr{}
%\email{}
%\thanks{}


%\keywords{rationally elliptic toric manifolds}

\date{\today}

%\dedicatory{}

\begin{abstract}
  In this note we classify rationally elliptic simply connected compact toric orbifolds up to algebraic isomorphism.
\end{abstract}

\maketitle

\section{Introduction}

In rational homotopy theory  it is shown that there are two types of simply connected spaces with finite dimensional rational cohomology: rationally elliptic and rationally hyperbolic spaces.
For rationally elliptic spaces \(X\) the total dimension \(\sum_{i\geq 2}\dim \pi_i(X)\otimes \mathbb{Q}\) of the rational homotopy groups is finite, whereas for rationally hyperbolic spaces the sum \(\sum_{i=2}^k \dim \pi_i(X)\otimes \mathbb{Q}\) grows exponentially (see for example \cite{MR1802847}).

A toric variety \(X\)
of complex dimension \(n\) is a normal complex algebraic variety with an
action of a complex torus \((\mathbb{C}^*)^n\) having an open dense orbit.
If \(X\) is compact and smooth we call it a toric manifold.

In the recent paper \cite{biswas19:_ration} rationally elliptic toric manifolds in complex dimension at most three were classified up to algebraic isomorphism.
In toric topology generalizations of toric varieties such as torus manifolds and torus orbifolds are studied.
A classification of rationally elliptic torus orbifolds up to rational homotopy equivalence has been given in \cite{MR3862119}.
The aim of this note is to explain how the methods of the latter paper lead to a classification result (up to algebraic isomorphism) for rationally elliptic toric manifolds and orbifolds in all dimensions.
Our main result is as follows.

\begin{theorem}
\label{sec:introduction}
  Let \(X\) be a compact simply connected toric orbifold of complex dimension \(n\geq 1\) which is rationally elliptic.
  Then there is an algebraic isomorphism \(X\rightarrow X'\) where \(X'\) is a quotient of an almost free action of an abelian complex algebraic group \(G\) on \(Y=\prod_i(\mathbb{C}^{n_i+1}-\{0\})\), for certain \(n_i>0\) (depending on \(X\)) with \(\sum_i n_i=n\).
\end{theorem}

In case that \(X\) is a toric manifold, \(G\) is a complex torus acting freely on \(Y\).
Therefore it follows that \(X'\) is a so-called generalized Bott manifold.

Generalized Bott manifolds are certain projective toric manifolds \cite[p. 300-302]{MR3363157}. They can be constructed as total spaces of towers of fiber bundles
\begin{equation*}
  X=X_n\rightarrow X_{n-1}\rightarrow\dots\rightarrow X_1\rightarrow X_0=\{pt\},
\end{equation*}
where each \(X_i\) is the projectivization of a Whitney sum of complex line bundles over \(X_{i-1}\).
Generalized Bott manifolds have been studied intensively by the Japanese--Korean school of toric topologists (see for example \cite{MR2666127}, \cite{MR2792373}, \cite{MR2983437}, \cite{MR3342682}, \cite{MR3391366}, \cite{MR3342675}, \cite{MR3273878}).

We note that all generalized Bott manifolds are rationally elliptic so that the above theorem gives a complete classification for rationally elliptic toric manifolds.

The proof of this result combines the quotient construction of toric varieties due to Cox \cite{MR1299003} with a recent result on the combinatorics of orbit spaces of rationally elliptic torus orbifolds given in \cite{MR3862119}.
Note that in the manifold case the arguments of \cite{MR3862119} also hold for torus manifolds with invariant metrics of non-negative sectional curvature (see \cite{MR3355120}).
Therefore Theorem~\ref{sec:introduction} also holds for toric manifolds admitting a non-negatively curved Riemannian metric which is invariant under the action of the maximal compact torus in \((\mathbb{C}^*)^n\).

This note has three more sections. In the next Section~\ref{sec:toric-vari-quas} we recall the construction of toric varieties as quotient spaces.
Then in Section~\ref{sec:torus_manifolds}, we recall the classification of rationally elliptic torus manifolds and orbifolds up to homeomorphism and rational homotopy equivalence given in \cite{MR3355120} and \cite{MR3862119}. 
In the last Section~\ref{sec:proof-theor-refs} we prove Theorem~\ref{sec:introduction}.

\section{The quotient construction of toric varieties}
\label{sec:toric-vari-quas}

In this section we recall the basic notions of toric geometry and describe the quotient construction of toric varieties.

For an introduction to toric geometry we refer the reader to \cite{MR2810322}, \cite{MR1234037} and \cite{MR922894}.

A toric variety \(X\)
of complex dimension \(n\) is a normal complex algebraic variety with an
action of a complex torus \((\mathbb{C}^*)^n\) having an open dense orbit.
If \(X\) is compact and smooth we call it a toric manifold.

The equivariant isomorphism types of these varieties are in one-to-one correspondence with combinatorial objects called fans (see for example \cite[Section 3.1]{MR2810322}).
A fan \(F\) is a finite collection of convex polyhedral cones in \(\mathbb{R}^n\) such that all faces of a cone \(C\in F\) are again in \(F\) and the intersection of any two cones \(C_1,C_2\in F\) is a face of each \(C_1\) and \(C_2\).

If \(X\) is a toric variety and \(F_X\) the fan corresponding to \(X\), then, by \cite[Theorem 3.2.6]{MR2810322}, there is a bijection \(O\mapsto C_O\) between the set of \((\mathbb{C}^*)^n\)-orbits in \(X\) and the set of cones in \(F_X\), such that
\begin{enumerate}
\item \(\codim_{\mathbb{C}} O=\dim_{\mathbb{R}} C_O\) for all orbits \(O\subset X\),
  \item if \(O_1\) and \(O_2\) are orbits in \(X\), then \(O_1\) is contained in the closure of \(O_2\) if and only if \(C_{O_2}\) is a face of \(C_{O_1}\).
\end{enumerate}

A \(k\)-dimensional cone \(C\) is called simplicial if it is spanned by \(k\) linearly independent vectors \(v_1,\dots,v_k\in \mathbb{R}^n\).
In case that the simplicial cone \(C\) belongs to a fan \(F\), the rays spanned by the \(v_i\) also belong to \(F\).
A toric variety \(X\) is an orbifold if and only if its corresponding fan \(F_X\) is simplicial, i.e. all its cones are simplicial.
\(X\) is compact if and only if the union of all cones in \(F_X\) is \(\mathbb{R}^n\) \cite[Theorem 3.1.19]{MR2810322}.

From a simplicial fan \(F_X\) we can construct an abstract simplicial complex \(\Sigma_X\)  as follows.
Let \(I\) be the set of rays of \(F_X\).
A subset \(\sigma=\{i_1,\dots i_k\}\subset I\) is a simplex of \(\Sigma_X\) if and only if the rays \(i_1,\dots,i_k\) span a \(k\)-dimensional cone in \(F_X\).
Note that there is a natural one-to-one correspondence between the cones of \(F_X\) of positive dimension and the simplices of \(\Sigma_X\).

Let \(T=(S^1)^n\subset (\mathbb{C}^*)^n\) be the maximal compact torus.
In case that \(X\) is a compact toric orbifold, the above simplicial complex can also be described in terms of the stratification of \(X/T\) by the identity components of the isotropy groups of the \(T\)-action on \(X\).

This goes as follows. For a closed connected subgroup \(H\subset T\), let \(S_H\) be the set of orbits \(Tx\in X/T\) of types \(T/T_x\) such that the identity component of \(T_x\) is equal to \(H\).
We call the connected components of \(S_H\) the open \(H\)-strata of \(X/T\).
The closed \(H\)-strata are the closures of the open \(H\)-strata.
The codimension of an (closed or open) \(H\)-stratum is the dimension of \(H\).

Since \(X\) is compact and the \(T\)-action has only finitely many orbit types, the set \(\mathcal{P}(X/T)\)  of all closed strata of positive codimension is finite.
Moreover, it is partially ordered by inclusion.
Therefore \(\mathcal{P}(X/T)\) is a poset, the so-called face poset of \(X/T\).

It is easy to see that the open codimension-\(k\) strata of \(X/T\) are given by the subsets \(O/T\subset X/T\), where \(O\) runs through the \((\mathbb{C}^*)^n\)-orbits of (complex) codimension \(k\) in \(X\).
Therefore \(\mathcal{P}(X/T)\) is dual to the simplicial complex \(\Sigma_X\) in the following sense:
There is an order reversing bijection \(\Sigma_X\rightarrow \mathcal{P}(X/T)\) such that the \((k-1)\)-dimensional simplices of \(\Sigma_X\) correspond to the codimension-\(k\) strata of \(X/T\).
Here the simplicial complex \(\Sigma_X\) is also partially ordered by inclusion of simplices.

In particular, the intersection \(S_1\cap S_2\)  of any two closed strata \(S_1,S_2\)  of \(X/T\) is connected or empty.
This follows from the fact that \(S_1\cap S_2\) is the disjoint union of those closed strata which are maximal among those strata which are contained in both \(S_1\) and \(S_2\). Note here that for any two simplices \(\sigma_1\) and \(\sigma_2\)  in \(\Sigma_X\) there is at most one minimal simplex which contains both \(\sigma_1\) and \(\sigma_2\).
If such a simplex exists in \(\Sigma_X\), then it is given by \(\sigma_1\cup\sigma_2\).

Cox \cite{MR1299003} (see also \cite[Chapter 5]{MR2810322}) gave a description of \(X\) as a quotient of an almost free action by an abelian complex algebraic group \(G\) on an open dense subset \(Y(\Sigma_X)\) of \(\mathbb{C}^I\).

The set \(Y(\Sigma_X)\) can be defined as follows.
For \(z\in \mathbb{C}^I\) let \(I(z)=\{i\in I;\; z(i)=0\}\). We then define
\begin{equation*}
  Y(\Sigma_X)=\{z\in \mathbb{C}^I;\;I(z)\in \Sigma_X\cup \{\emptyset\}\}.
\end{equation*}

With this notation Cox's description of a toric orbifold as a quotient can be stated as follows.

\begin{theorem}
  \label{sec:quot-constr-toric}
  Let \(X\) be a toric orbifold.
  Then \(X\) is algebraically isomorphic to a quotient of an almost free action of an abelian complex algebraic group \(G\) on \(Y(\Sigma_X)\).
Moreover, in case \(X\) is a toric manifold then \(G\) is a complex torus which acts freely on \(Y(\Sigma_X)\).
\end{theorem}

We close this section by giving two examples of sets \(Y(\Sigma_X)\) for special choices of compact toric orbifolds \(X\).

\begin{example}
  \label{sec:quot-constr-toric-1}
  Let \(X\) be a compact toric orbifold.
  If \(\Sigma_X\) is dual to the face poset of an \(n\)-dimensional simplex \(\Delta^n\), then \(Y(\Sigma_X)=\mathbb{C}^{n+1}-\{0\}\). This is for example the case if \(X=\mathbb{C} P^n\).
\end{example}

\begin{example}
  \label{sec:quot-constr-toric-2}
  Let \(X_1,X_2,X_3\) be compact toric orbifolds of complex dimensions \(n_1,n_2,n_3\), respectively, such that \(n_1=n_2+n_3\).
  We equip \(X_2\times X_3\) with the natural product action by \((\mathbb{C}^*)^{n_1}=(\mathbb{C}^*)^{n_2}\times (\mathbb{C}^*)^{n_3}\).
  With this action \(X_2\times X_3\) becomes a compact toric orbifold of complex dimension \(n_1=n_2+n_3\).
  
  If \(\mathcal{P}(X_1/T^{n_1})\) is isomorphic to \(\mathcal{P}((X_2\times X_3)/T^{n_1})\) then \(Y(\Sigma_{X_1})\cong Y(\Sigma_{X_2})\times Y(\Sigma_{X_3})\).
  Note also that \((X_2\times X_3)/T^{n_1}\) is strata preserving homeomorphic to \((X_2/T^{n_2})\times (X_3/T^{n_3})\).
\end{example}

\section{Rationally elliptic torus manifolds revisited}
\label{sec:torus_manifolds}

In this section we recall the definition of torus manifolds and orbifolds and classification results for simply connected rationally elliptic torus manifolds and orbifolds.
Torus manifolds and orbifolds are studied in toric topology.
Toric topology has its origin in the paper \cite{MR1104531}. We refer the reader to \cite{MR1897064} and \cite{MR3363157} for an overview over the development of the subject since then.

A torus manifold \(M\) is a closed, connected, orientable manifold of (real) dimension \(2n\) equipped with an effective action of an \(n\)-dimensional compact torus \(T=(S^1)^n\), such that the fixed point set \(M^{T}\) is non-empty.
Torus orbifolds are  natural generalizations of torus manifolds. One gets their definition if one replaces the word ``manifold'' by the word ``orbifold'' in the above definition of a torus manifold.

Note that toric manifolds (and compact toric orbifolds) equipped with the action of the maximal compact torus \(T=(S^1)^n\subset (\mathbb{C}^*)^n\) are torus manifolds (and torus orbifolds, respectively).
Note, moreover, that the definition of the face poset of the orbit space of a compact toric orbifold carries over without changes to the situation of a torus orbifold.
However, in this case it is no longer true that the face poset is always dual to a simplicial complex.

A torus manifold \(M\) is called locally standard if the \(T\)-action on \(M\) is locally modeled on effective \(T\)-representations on \(\mathbb{C}^n\).
In this case the orbit space \(M/T\) is naturally a nice manifold with corners and all isotropy groups are connected.

Here a manifold with corners is called nice, if all its codimension-\(k\) faces are contained in exactly \(k\) codimension-one faces.
In this case each codimension-\(k\) face is a component of the intersection of exactly \(k\) codimenion one faces (see \cite{MR2283418} or \cite{MR3030690} for more details).

The stratification of \(M/T\) by orbit types coincides with the stratification of \(M/T\) by faces.
In particular, the poset \(\mathcal{P}(M/T)\) is the poset of faces of \(M/T\) (viewed as a nice manifold with corners).
This justifies the name face poset for \(\mathcal{P}(M/T)\).

In \cite{MR3355120} simply connected rationally elliptic torus manifolds \(M\) with \(H^{\text{odd}}(M;\mathbb{Z})=0\) have been classified up to homeomorphism. These manifolds are all homeomorphic to quotients of free torus actions on products of spheres.

Note that by \cite{MR2283418} the cohomological condition implies that \(M\) is locally standard and that the faces \(F\) of \(M/T\) are acyclic over the integers, i.e. \(\tilde{H}^*(F;\mathbb{Z})=0\) for all faces \(F\) of \(M/T\).

The proof of the classification result in \cite{MR3355120} proceeds in two steps.
First it has been shown that the homeomorphism type of a torus manifold \(M\) as above depends only on combinatorial data, namely on the face poset \(\mathcal{P}(M/T)\) and the characteristic function \(\lambda\), which assigns to a face of \(M/T\) the isotropy group of a generic orbit in that face \cite[Theorem 3.4]{MR3355120}.

In a second step it has been shown, by a combinatorial argument (see Proposition 4.5 in \cite{MR3355120} and Theorem~\ref{sec:rati-ellipt-torus} below), that \(\mathcal{P}(M/T)\) is isomorphic to the face poset of a product
\[\prod_{i<r} \Sigma^{n_i}\times \prod_{i\geq r} \Delta^{n_i},\]
where \(\Delta^k\) is a \(k\)-dimensional simplex and \(\Sigma^k\) is the suspension of \(\Delta^{k-1}\).
In particular, \(\Delta^k\) and \(\Sigma^k\) are orbit spaces of the natural action of a maximal torus of the orthogonal group \(O(2k+2)\) (\(O(2k+1)\), respectively) on spheres of dimensions \(2k+1\) and \(2k\), respectively.

Note that \(\Delta^k\) has \(k+1\) codimension-one faces and the intersection of any \(k'\) of them is non-empty and connected for any \(k'\leq k\) and empty for \(k'=k+1\).
Therefore \(\Sigma^k\) has exactly \(k\) codimension-one faces and the intersection of any \(k'\) of them is connected for \(k'<k\) and contains exactly two isolated points if \(k'=k\).

By combining the above two steps, it follows that \(M\) is homeomorphic to a locally standard torus manifold \(M'\) with \(M'/T=\prod_{i<r} \Sigma^{n_i}\times \prod_{i\geq r} \Delta^{n_i}\).
Since each such \(M'\) is the quotient of a free torus action on a product of spheres, it follows that \(M\) is homeomorphic to such a quotient.

In \cite{MR3862119} this argument was generalized to simply connected rationally elliptic torus manifolds \(M\) with \(H^{\text{odd}}(M;\mathbb{Z})\neq 0\) and to simply connected rationally elliptic torus orbifolds \(\mathcal{O}\).
  While the combinatorial part of the proof goes through with modifications, the proof of the first step does not.
  Therefore in \cite{MR3862119} we only get a classification up to rational homotopy equivalence. Indeed all such \(\mathcal{O}\) are rationally homotopy equivalent to quotients of almost free torus actions on products of spheres.
  
  However, from the discussion in that paper we have the following combinatorial result.
  
\begin{theorem}
  \label{sec:rati-ellipt-torus}
    Let \(\mathcal{O}^{2n}/T\) be the orbit space of a simply connected rationally elliptic torus orbifold \(\mathcal{O}^{2n}\). Then \(\mathcal{P}(\mathcal{O}^{2n}/T)\) is isomorphic to the face poset of a product \(\prod_{i< r}\Sigma^{n_i}\times \prod_{i\geq r} \Delta^{n_i}\) with \(n_i>0\) and \(\sum_i n_i=n\).
\end{theorem}

The original proof of this result, given in \cite{MR3355120} and \cite{MR3862119}, was very long and technical.
A much simpler proof has later been given in \cite[Section 8]{goertsches18:_non_gkm}.

\section{The proof of Theorem \ref{sec:introduction}}
\label{sec:proof-theor-refs}

Now assume that \(X\) is a simply connected rationally elliptic compact toric orbifold of complex dimension \(n\).
Then \(X\) is, in particular, a  \(2n\)-dimensional torus orbifold.
Therefore by Theorem~\ref{sec:rati-ellipt-torus}, \(\mathcal{P}(X/T)\) is isomorphic to the face poset of a product
\[\prod_{i<r} \Sigma^{n_i}\times \prod_{i\geq r} \Delta^{n_i},\]
with \(n_i>0\) and \(\sum_i n_i=n\).

Note that the intersection of any two one-dimensional faces in \(\Sigma^{k}\), \(k>1\), is disconnected.
Therefore, (and because \(\Sigma^1=\Delta^1\)) all factors in the above product are of type \(\Delta^k\).

In other words, \(\Sigma_X\) is dual to the face poset of a product of simplices.
In particular, by Examples~\ref{sec:quot-constr-toric-1} and \ref{sec:quot-constr-toric-2}, we have
\[Y(\Sigma_X)=\prod_i(\mathbb{C}^{n_i+1}-\{0\}).\]
Hence, our Theorem~\ref{sec:introduction} follows from Theorem~\ref{sec:quot-constr-toric}.

\end{document}